\documentclass[a4paper,11pt,twoside]{article}
\usepackage{ntheorem}
\theoremstyle{break}
\usepackage{titling}
\newtheorem{theorem}{Theorem}

\newtheorem{definition}[theorem]{Definition}
\newtheorem{properties}[theorem]{Property}

\newenvironment{proof}[1][Proof]{\textbf{#1.} }{\hfill $\square$}

\usepackage{tikz}
\usepackage{latexsym}
\usepackage{amsmath}
\usepackage{amssymb}
\usepackage{hyperref}
\hypersetup{colorlinks=true,linkcolor=blue,urlcolor=blue,citecolor=blue}
\usepackage{titlesec}
\usepackage{cite}
\usepackage{enumitem}
\usepackage{url}
\usepackage{vmargin}
\setmarginsrb{2cm}{1.5cm}{2cm}{2.3cm}{0mm}{5mm}{0mm}{10mm}
\usepackage[utf8]{inputenc}
\titleformat*{\section}{\normalsize \bfseries}
\titleformat*{\subsection}{\small \bfseries}
\usepackage[font=small,labelfont=small]{caption}
\usepackage{fancyhdr}

\fancypagestyle{firststyle}
{
\fancyhead[LE,RO]{}
\fancyhead[RE,LO]{\textsl{\small Publication Scientifique de l'AEIF \ 2(2002) 127--132}\\ \textsl{\small \copyright \ PubScieAEIF}}
\fancyfoot{}
}

\def\thebibliography#1{\setlength{\itemsep}{0pt}\setlength{\parsep}{0pt}
  \section*{
  {}{REFERENCES}}\list 
  {[\arabic{enumi}]}{\settowidth\labelwidth{[#1]}\leftmargin\labelwidth
    \advance\leftmargin\labelsep
    \usecounter{enumi}}
    \def\newblock{\hskip .11em plus .33em minus .07em}
    \sloppy\clubpenalty4000\widowpenalty4000
    \sfcode`\.=1000\relax}

\begin{document}
\title{\bf{\large PROPERTIES OF CHEBYSHEV POLYNOMIALS}}
\author{{\normalsize Natanael Karjanto}\\ 
{\small \itshape Department of Applied Mathematics, University of Twente}\\
{\small \itshape P.O. Box 217, 7500 AE Enschede, The Netherlands}\\
{\small \protect \url{n.karjanto@math.utwente.nl}}}
\date{}
\maketitle

\thispagestyle{firststyle}

\begin{abstract}
\noindent
Ordinary differential equations and boundary value problems arise in many aspects of mathematical physics. Chebyshev differential equation is one special case of the Sturm-Liouville boundary value problem. Generating function, recursive formula, orthogonality, and Parseval's identity are some important properties of Chebyshev polynomials. Compared with a Fourier series, an interpolation function using Chebyshev polynomials is more accurate in approximating polynomial functions.\\

\end{abstract}

\def\abstractname{\textit{\textbf{R\'{e}sum\'{e}}}}
\begin{abstract}
\noindent	
\emph{Des équations différentielles ordinaires et des problèmes de valeurs limites se posent dans de nombreux aspects de la physique mathématique. L'équation différentielle de Chebychev est un cas particulier du problème de la valeur limite de Sturm-Liouville. La fonction génératrice, la formule récursive, l'orthogonalité et l'identité de Parseval sont quelques propriétés importantes du polynôme de Chebyshev. Par rapport à une série de Fourier, une fonction d'interpolation utilisant des polynômes de Chebyshev est plus précise dans l'approximation des fonctions polynomiales.}\\

\noindent
\emph{{\bfseries Mots-clés}: problème de valeur limite de Sturm-Liouville, équation différentielle de Chebyshev, polynôme de Chebyshev, fonction génératrice, formule récursive, orthogonalité, identité de Parseval.}
\end{abstract}

\section{INTRODUCTION}

Based on its name, these polynomials were first investigated by a Russian mathematician Pafnuty Lvovich Chebyshev (1821-1894). In addition to orthogonal functions, Chebyshev also studied inequalities, prime numbers, probability theory, quadratic forms, integral theory, geographical map, geometrical volume formulas, and mechanics as well as problems in changing circular motion to straight motion with mechanical coupling~\cite{britannica86}. In this paper, we will observe the appearance of Chebyshev polynomials as solutions of one particular case from a Sturm-Liouville boundary value problem. We discuss and prove several essential properties, such as the generating formula, the recursive relation, and Parseval's identity. We cover but do not prove Rodrigues' formula, even and odd functions, as well as orthogonality. Depending on the type of functions, an interpolation function using Chebyshev polynomials has a faster rate of convergence and better accuracy than a Fourier series when approximating a polynomial. On the other hand, a Fourier series converges faster and better when interpolating non-polynomial functions.

\section{CHEBYSHEV DIFFERENTIAL EQUATION AND ITS POLYNOMIALS}
\begin{definition}[Sturm-Liouville boundary value problem]
An ordinary differential equation (ODE) defined on an interval $a \leq x \leq b$ with a general form
\begin{equation}
\frac{d}{dx}\left(p(x)\,\frac{dy}{dx}\right) + \left[q(x) + \lambda \,r(x) \right] \,y = 0 			\label{STBVP}
\end{equation}
and boundary conditions
\begin{equation}
\left\{ 
\begin{aligned}
a_{1}\,y(a) + a_{2}\,y'(a) &= 0  \\
b_{1}\,y(b) + b_{2}\,y'(b) &= 0
\end{aligned} 
\right.
\end{equation}
is called a {\upshape Sturm-Liouville boundary value problem} or a {\upshape Sturm-Liouville system} with $p(x) > 0$, $q(x)$, and a weighted function $r(x) > 0$ are given functions; $a_{1}$, $a_{2}$, $b_{1}$, and $b_{2}$ are given constants; and the characteristic value $\lambda$ is an undetermined parameter.
\end{definition}

The boundary value problem~\eqref{STBVP} admits several special cases that depend on the defined interval, given functions, and characteristic values~\cite{efunda02}. Each of these special cases leads to a particular ODE with orthogonal functions as its solutions. Among them are Bessel's functions, Legendre polynomials, Hermite polynomials, Laguerre polynomials, and Chebyshev polynomials. This paper discusses the latter in particular.

For the values of an endpoint interval $a = -1$ and $b = 1$, the given functions $p(x) = \sqrt{1 - x^2}$, $q(x) = 0$, $r(x) = 1/\sqrt{1-x^2}$, and the characteristic value $\lambda = n^2$, the ODE~\eqref{STBVP} becomes
\begin{equation}
\frac{d}{dx}\left(\sqrt{1 - x^2}\,\frac{dy}{dx}\right) + \frac{n^2}{\sqrt{1 - x^2}}\,y = 0. 		\label{ChebyBVP}
\end{equation}

\begin{definition}[Chebyshev differential equation]
The ODE~\eqref{ChebyBVP}, or expressed in another form as
\begin{equation}
(1 - x^2)\frac{d^2y}{dx^2} - x\frac{dy}{dx} + n^2\,y = 0			\label{PDCheby}
\end{equation}
with $n \in \mathbb{N}_0$ for $|x| < 1$ is called the {\upshape Chebyshev differential equation}. The solution for this ODE is called the {\upshape Chebyshev function} with singularities at $x = -1, 1$, and $\infty$.
\end{definition}
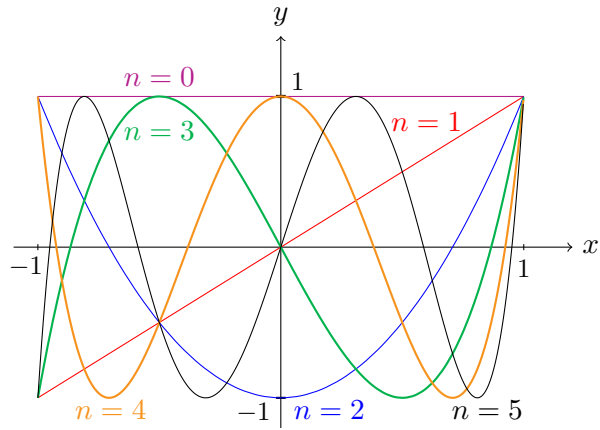
\begin{figure}[htbp]
\vspace*{-0.2cm}
\begin{center}
\begin{tikzpicture}[domain=-1:1, yscale=2.5, xscale = 4, scale = 0.8]
\draw[->] (-1.1,0) -- (1.2,0) node[right] {$x$};
\draw[->] (0,-1.2) -- (0,1.4) node[above] {$y$};
\draw (1,0.02) -- (1,-0.02) node[below] {\small $1$};
\draw (-1,0.02) -- (-1,-0.02);
\draw (-1.05,0) node[below] {\small $-1$};
\draw[very thick] (0.02,1) -- (-0.02,1);
\draw (0,1.1) node[right] {\small 1};
\draw[thick] (0.02,-1) -- (-0.02,-1);
\draw (0,-1.1) node[left] {\small $-1$};
\draw[smooth,variable=\x,magenta!70!blue] plot ({\x},{1});
\draw[magenta!70!blue] (-0.5,1) node[above] {$n = 0$};
\draw[smooth,variable=\x,red] plot ({\x},{\x});
\draw[red] (0.6,0.7) node[above] {$n = 1$};
\draw[smooth,variable=\x,blue] plot ({\x},{2*\x*\x - 1});
\draw[blue] (0.2,-0.95) node[below] {$n = 2$};
\draw[smooth,thick,variable=\x,green!70!blue]  plot[samples = 50] ({\x},{4*\x*\x*\x - 3*\x});
\draw[green!70!blue] (-0.5,0.9) node[below] {$n = 3$};
\draw[smooth,thick,variable=\x,yellow!50!red] plot[samples = 100] ({\x},{8*\x*\x*\x*\x - 8*\x*\x + 1});
\draw[yellow!50!red] (-0.7,-0.95) node[below] {$n = 4$};
\draw[smooth,variable=\x,black] plot[samples = 200] ({\x},{16*\x*\x*\x*\x*\x - 20*\x*\x*\x + 5*\x});
\draw[black] (0.85,-0.95) node[below] {$n = 5$};
\end{tikzpicture}
\end{center}
\caption{Plot of Chebyshev polynomials $T_n(x)$ for $n = 0, 1, 2, \dots, 5$.}  \label{plot}
\end{figure}
\begin{definition}[Chebyshev polynomials]
The functions
\begin{equation}
T_{n}(x) = \cos\,(n\,\cos^{-1} x) \qquad \text{\upshape or} \qquad
T_{n}\,(\cos \theta) = \cos\,(n\,\theta) \qquad \text{\upshape for}	\qquad x = \cos \theta
\end{equation}
are called the {\upshape Chebyshev polynomials} or the {\upshape Chebyshev polynomials of the first kind}~\cite{spiegel74,andrews92,mason02}.
These are solutions of the ODE~\eqref{PDCheby} for non-negative integers $n$, $n \in \mathbb{N}_0$. See Figure~\ref{plot}.
\end{definition}
Using a series expansion ${\displaystyle y(x) = \sum_{n = 0}^{\infty} \, a_{n} \,x^n}$, a general solution for the ODE~\eqref{PDCheby} is givey by
\begin{equation}
y(x) = b_{1} \, T_{n}(x) + b_{2} \, \sqrt{1 - x^2}\,U_{n - 1}(x)
\end{equation}
where $U_n(x) = \sin\left[(n + 1) \cos^{-1} x \right]/\sin(\cos^{-1} x)$ are the Chebyshev polynomials of the second kind.

\section{PROPERTIES OF THE CHEBYSHEV POLYNOMIALS}
\begin{properties}[Rodrigues' formula]
The Chebyshev polynomials $T_n(x)$ can be expressed as Rodrigues' formula in terms of derivatives	
\begin{equation}
T_{n}(x) = \frac{\sqrt{1 - x^2}}{(-1)^n\,(2n-1)(2n-3)\dots
1}\frac{d^n}{dx^n}(1-x^2)^{n-\frac{1}{2}}
\end{equation}
with $n = 0, 1, 2,  3, \dots$.
\end{properties}

\begin{properties}[Generating function] 		 \label{generating}
The generating function for the Chebyshev polynomials $T_n(x)$ reads
\begin{equation}
\frac{1 - z\,x}{1 - 2\,z\,x + z^2} = \sum_{n = 0}^{\infty} \, T_{n}(x) \, z^n.
\end{equation}
\end{properties}
\begin{proof}
By letting $x = \cos \theta$, we have
\begin{align*}
\sum_{n=0}^{\infty}\,T_{n}(x)\,z^n &=
\sum_{n=0}^{\infty}\,T_{n}(\cos \theta)\,z^n =
\sum_{n=0}^{\infty}\,\cos (n\,\theta) \,z^n = 1 + \frac{1}{2}\,
\sum_{n = -\infty,\,n \neq 0}^{\infty}\,z^{|n|}\,e^{\,i\,n\,\theta} \\ &= 
1 + \frac{1}{2}\,\left( \sum_{n=1}^{\infty}\,z^n\,e^{\,i\,n\,\theta} +
\sum_{n=1}^{\infty}\,z^n\,e^{-i\,n\,\theta} \right) = 1 +
\frac{1}{2}\,\left(\frac{z\,e^{\,i\,\theta}}{1 -
z\,e^{\,i\,\theta}} + \frac{z\,e^{-i\,\theta}}{1 - z\,e^{-i\,\theta}} \right) \\ &= 
1 + \frac{1}{2}\,\left(\frac{z\,e^{\,i\,\theta} - z^2 + z\,e^{-i\,\theta} -z^2}{1 - z\,e^{\,i\,\theta}- z\,e^{-i\,\theta} + z^2} \right) = 
1 + \frac{z\,\cos \theta - z^2}{1 - 2\,z\,\cos \theta + z^2} \nonumber \\ &= 
1 + \frac{x\,z}{1 - 2\,x\,z + z^2} = \frac{1 - x\,z}{1 - 2\,x\,z + z^2}.
\end{align*}
Hence, we have shown the generating function formula for the Chebyshev polynomials.
\end{proof}

\begin{properties}[Even and odd functions]
Based on their degrees, the Chebyshev polynomials $T_{n}(x)$ admit two possibilities as either even or odd functions.
\begin{itemize}[leftmargin=1.2em]
\item For $n$ even, then $T_{n}(x)$ are even functions.

\item For $n$  odd, then $T_{n}(x)$ are  odd functions.
\end{itemize}
\end{properties}

\begin{properties}[Recursive relationship]
One Chebyshev polynomial $T_{n}(x)$ at a particular point can be expressed in terms of its neighboring Chebyshev polynomials at the identical point
\begin{equation}
T_{n + 1}(x) = 2 \,x\,T_{n}(x) - T_{n - 1}(x).			\label{recu}
\end{equation}
\end{properties}
We refer~\eqref{recu} as a three-term recursive relationship since the formula forms a relationship among three terms of successive Chebyshev polynomial. \\
\begin{proof}
From the definition of the Chebyshev polynomials, we have $T_{n}(\cos \theta) = \cos(n \theta)$. 
We also have
\begin{equation*}
T_{n+1}(x) = T_{n+1}(\cos \theta) = \cos (n+1)\,\theta = \cos n\theta \,\cos \theta - \sin n\theta \, \sin \theta
\end{equation*}
and 
\begin{equation*}
T_{n-1}(x) = T_{n-1}(\cos \theta) = \cos (n-1)\,\theta = \cos n\theta \,\cos \theta + \sin n\theta \, \sin \theta.
\end{equation*}
Next, by adding these two expressions, we obtain 
\begin{equation*}
T_{n+1}(x) + T_{n-1}(x) = 2\,\cos n\theta \,\cos \theta = 2\,T_{n}(x)\,x.
\end{equation*}
Finally, we attain the recursive formula $T_{n+1}(x) = 2\,x\,T_{n}(x) - T_{n-1}(x)$.
\end{proof} \\
This recursive relationship can also be derived from the generating function for the Chebyshev polynomials (Property~\ref{generating}) by differentiating it twice with respect to $z$ and equating the coefficients for $z^n$ on both sides of the equation. We will obtain the desired result accordingly.

\begin{properties}[Orthogonality] 			\label{ortogonalitas}
The Chebyshev polynomials $T_{n}(x)$ form a complete orthogonal set on the interval $-1 \leq x \leq 1$ with a weighted function $1/\sqrt{1 - x^2}$.
\end{properties}
This means
\begin{equation}
\int_{-1}^{1}\,\frac{1}{\sqrt{1-x^2}}\,T_{m}(x)\,T_{n}(x)\,dx =
\left\{ {\begin{array}{*{20}c} 0, \\ \pi, \\ \frac{\pi}{2}, \\
\end{array}} \right. \qquad
\left. {\begin{array}{l} m \neq n  \\ m = n = 0 \\ m = n = 1,
2, 3, \dots
\end{array}} \right.
\end{equation}
Using this orthogonality property, a piecewise continuous function on interval $-1 \leq x \leq 1$ can be expressed in terms of the Chebyshev polynomials
\begin{equation}
\sum_{n = 0}^{\infty} \, C_{n} \, T_{n}(x) = \left\{
{\begin{array}{*{20}c} f(x), \\ \frac{1}{2}\,[f(x^-) + f(x^+)], \\
\end{array}} \right. \qquad 
\left. {\begin{array}{l} \textmd{for continuous $f(x)$ } \\
\textmd{at discontinuous points}
\end{array}} \right.
\end{equation}
where
\begin{equation}
C_{0} = \frac{1}{\pi}\,\int_{-1}^{1}\,\frac{f(x)}{\sqrt{1 - x^2}} \, dx \qquad \textmd{and} \qquad
C_{n} = \frac{2}{\pi}\,\int_{-1}^{1}\,\frac{1}{\sqrt{1 - x^2}}\,f(x)\,T_{n}(x)\,dx, \qquad \textmd{for} \quad n \in \mathbb{N}.  \label{coeff}
\end{equation}

\begin{properties}[Parseval's identity] %\hfill
Parseval's identity related to the series expansion ${\displaystyle f(x) = \sum_{n=0}^{\infty}\,C_{n}\,T_{n}(x)}$ with the Chebyshev polynomials $T_{n}(x)$ is given by
\begin{equation}
\int_{-\infty}^{\infty}\,\frac{1}{\sqrt{1 -x^2}}\,[f(x)]^2\,dx = C_{0}^2\,\pi + \frac{\pi}{2}\,\sum_{n=1}^{\infty}\,C_{n}^2.
\end{equation}
\end{properties}
\begin{proof}
The first step is by taking the square of both sides of the series expansion ${\displaystyle f(x) =
\sum_{n=0}^{\infty}\,C_{n}\,T_{n}(x)}$ to yield
\begin{equation}
[f(x)]^2 = \sum_{n=0}^{\infty}\,\sum_{m=0}^{\infty}\, C_{n}\, C_{m}\,T_{n}(x)\,T_{m}(x).   \label{squared}
\end{equation}
We then multiply~\eqref{squared} with the weighted function $1/\sqrt{1 - x^2}$ and integrate it with respect to the variable $x$ from $-\infty$ to $\infty$ to obtain
\begin{equation}
\int_{-\infty}^{\infty}\,\frac{1}{\sqrt{1 - x^2}}\,[f(x)]^2\,dx = 
\sum_{n=0}^{\infty}\,\sum_{m=0}^{\infty}\, C_{n}\, C_{m}\,\int_{-\infty}^{\infty}\,\frac{1}{\sqrt{1 - x^2}}T_{n}(x)\,T_{m}(x)\,dx. \label{integrated}
\end{equation}
Using Property~\ref{ortogonalitas}, we can write~\eqref{integrated} as follows:
\[
\int_{-\infty}^{\infty}\,\frac{1}{\sqrt{1 - x^2}}\,[f(x)]^2\,dx =
C_{0}^2\,\pi + \frac{\pi}{2}\,\sum_{n=1}^{\infty}\,C_{n}^2
\]
which gives Parseval's identity for the Chebyshev polynomials.
\end{proof}

\section{APPLICATION}

The Chebyshev polynomials play an important role in function approximation, integral approximation, and polynomial interpolation problem. For a periodic function $f(x)$ on a particular interval, a Fourier series is more superior in approximating it with a relatively tiny error. However, in applications, we often encounter non-periodic functions which are only defined in a limited interval. The Chebyshev polynomials have a substantial contribution to these cases~\cite{mason02,thompson94}. By approximating the function $f(x)$ with the Chebyshev polynomials $T_n(x)$, we obtain the following interpolation function $P_N(x)$:
\begin{equation}
f(x) \approx P_{N}(x) = \sum_{n = 0}^{N} \, C_{n} \, T_{n}(x).
\end{equation}
We must select the coefficients $C_n$ such that the norm $\|f(x) - P_{n}(x)\|_{2}$ is minimum.
By seeking extreme values (in this case minimum values) of the $L^2$-norm, we obtain the coefficients $C_n$ as in~\eqref{coeff}.
\begin{figure}[h]
\vspace*{-0.5cm}
\centering
\includegraphics[width=0.45\textwidth]{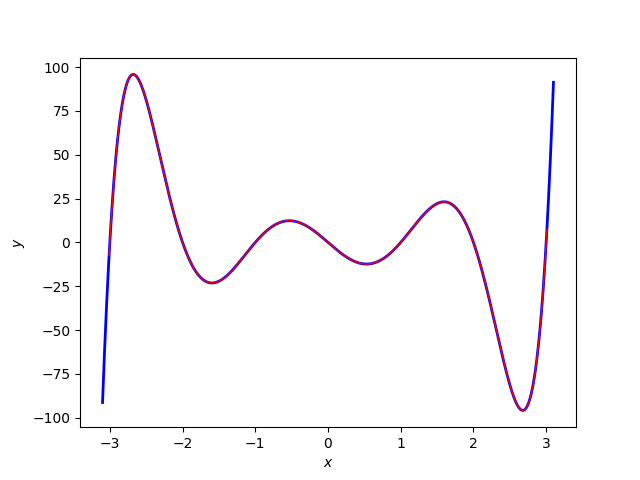} \hspace{1cm}
\includegraphics[width=0.45\textwidth]{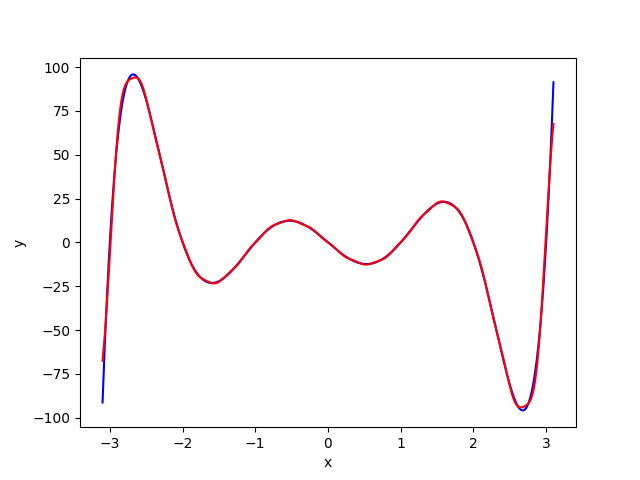}
\caption{Plot of a seventh-order polynomial $y = f(x) = x^7 - 14x^5 + 49x^3 - 36x$ (blue curve) approximated by interpolation functions corresponding to the Chebyshev polynomials with $N = 7$ (red curve on the left panel) and a Fourier series with $N = 20$ (red curve on the right panel). We observe that for an $N^\text{th}$-order polynomial, a Chebyshev approximation with $N$ terms provides an accurate interpolation while more terms might be needed in a Fourier series approximation for a comparable qualitative accuracy.}  \label{poly}
\end{figure}

An $N^\text{th}$-order interpolation function using the Chebyshev polynomials provides an accurate approximation for any $N^\text{th}$-order polynomial with zero error. Approximating a polynomial using the Fourier series requires more terms to reach a qualitative comparable accuracy. Figure~\ref{poly} illustrates an example for a seventh-order polynomial $y = f(x) = x^7 - 14x^5 + 49x^3 - 36x$ when it is approximated by the Chebyshev polynomials and Fourier series. For the former, taking $N = 7$ in the interpolation function $P_N(x)$ provides an accurate interpolation. For the latter, the polynomial could be interpolated with a relatively tiny but visible error by taking at least $N = 10$.

\begin{figure}[h]
\centering
\includegraphics[width=0.45\textwidth]{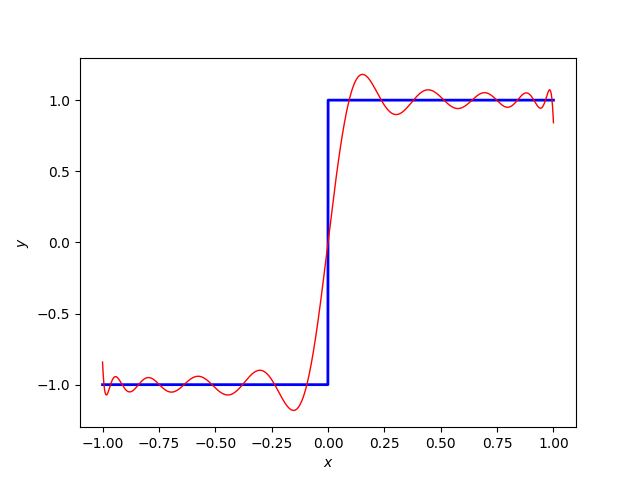} \hspace{1cm}
\includegraphics[width=0.45\textwidth]{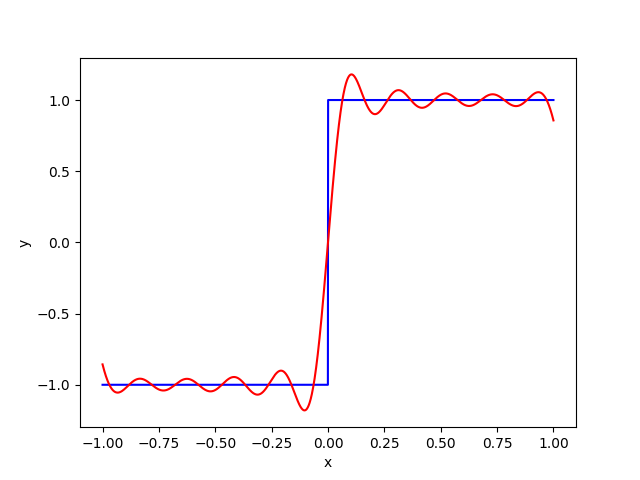}
\vspace*{-0.1cm}
\caption{A step function (blue curve) approximated by interpolation functions corresponding to the Chebyshev polynomials with $N = 20$ (red curve on the left panel) and a Fourier series with $N = 10$ (red curve on the right panel). For non-polynomial functions, more terms are required in the interpolation function for the Chebyshev polynomials to reach a similar qualitative accuracy with the one using the Fourier series.}  \label{step}
\end{figure}

For non-polynomial functions, we observe that an interpolation function using Fourier series converges faster than the one employing the Chebyshev polynomials. Figure~\ref{step} illustrates another example of a step (Heaviside) function approximated by the Chebyshev polynomials and Fourier series on the interval $-1 \leq x \leq 1$. We need up to $N = 20$ in the partial sum for the former while it is sufficient to include up to $N = 10$ for the latter to reach a qualitative comparable accuracy. Note that for this particular case, the chosen step function is odd, and so is the integrands in~\eqref{coeff} for $n$ is even. Thus, the coefficients $C_n = 0$ for $n$ is even. This means that although we take $N = 20$ for the interpolation function, only 10 nonzero terms contribute to the approximation, i.e., only the odd terms of the Chebyshev polynomial. For both cases, we observe the occurrence of the Gibbs phenomenon.\footnote[1]{The Gibbs phenomenon is an event where a jump or a deviation occurs near discontinous points of a piecewise continuous function when the function is approximated by Fourier series or other function series.}

\section{CONCLUSION}
In this paper, we have observed that the Chebyshev differential equation is a special case of the Sturm-Liouville boundary value problem. Its solutions, the Chebyshev polynomials, have properties which are also found in other polynomials. One of these properties is orthogonoality, which plays an important role in approximating functions, similar to the well-known Fourier series. Similar to the case in the Fourier series, approximating piecewise continuous functions using the Chebyshev polynomials also produces the Gibbs phenomenon. When approximating polynomials, the Chebyshev polynomials exhibit a superiority in comparison to the Fourier series due to its accuracy. However, the Fourier series seems to converge faster in interpolating non-polynomial functions than the Chebyshev polynomials.

\subsection*{ACKNOWLEDGMENT}
{\small The author thanked Hendra Gunawan (Bandung Institute of Technology), Intan Detiena Muchtadi-Alamsyah (The University of Picardy Jules Verne, Amiens), and Natanael Dewobroto (Paul Verlaine University, The University of Lorraine, Metz) for fruitful discussion, constructive advice, and valuable suggestion. \par}

\end{document}